\documentclass[12pt]{article}
\usepackage[utf8]{inputenc}
\usepackage{mathtools}
\usepackage[margin=1.15in]{geometry}
\usepackage{bm}
\usepackage{amsmath}
\usepackage{caption}
\usepackage{amsfonts}
\usepackage[svgnames]{xcolor}
\usepackage{float}
\usepackage{siunitx}
\usepackage{amsthm}
\usepackage{listings}
\usepackage{subfigure}
\usepackage{color}

\DeclarePairedDelimiter\floor{\lfloor}{\rfloor}
\pdfpageattr{/Group <</S /Transparency /I true /CS /DeviceRGB>>}
\theoremstyle{definition}
\newtheorem{definition}{Definition}[section]
\newtheorem{theorem}{Theorem}[section]
\definecolor{codegreen}{rgb}{0,0.6,0}
\definecolor{codegray}{rgb}{0.5,0.5,0.5}
\definecolor{codepurple}{rgb}{0.58,0,0.82}
\definecolor{backcolour}{rgb}{0.95,0.95,0.92}

\newcommand{\bmat}{\begin{matrix}}
\newcommand{\emat}{\end{matrix}}
\newcommand{\ben}{\begin{equation*}}
\newcommand{\een}{\end{equation*}}
\newcommand{\bean}{\begin{eqnarray*}}
\newcommand{\eean}{\end{eqnarray*}}
\newcommand{\be}{\begin{equation}}
\newcommand{\ee}{\end{equation}}
\newcommand{\bea}{\begin{eqnarray}}
\newcommand{\eea}{\end{eqnarray}}

\newcommand{\mbs}[1]{\boldsymbol{#1}}

\newcommand{\bx}{\mbs{x}}

\newcommand{\bp}{\mbs{p}}

\newcommand{\mcal}[1]{\mathcal{#1}}

\newcommand{\bpi}{{\mbs{\pi}}}
\newcommand{\lrp}[1]{\left(#1\right)}
\newcommand{\lrb}[1]{\left\{#1\right\}}

\newcommand{\given}{\,|\,}

\title{On the Uniformity of $(3/2)^n$ Modulo 1}
\author{Paula Neeley, Daniel Taylor-Rodriguez, J.J.P. Veerman, Thomas Roth}
\date{}
\begin{document}
\maketitle

\begin{abstract}
It has been conjectured that the sequence $(3/2)^n$ modulo 1 is uniformly distributed. The distribution of this sequence is significant in relation to unsolved problems in number theory including the Collatz conjecture. In this paper, we describe an algorithm to compute $(3/2)^n$ modulo 1 to $n=10^8$. We then statistically analyze its distribution. Our results strongly agree with the hypothesis that $(3/2)^n$ modulo 1 is uniformly distributed.
\end{abstract}


\section{Introduction}
\label{Introduction}
The distribution of $(3/2)^n$ modulo 1 is an interesting topic because of its connections with Mahler's Z-numbers \cite{mahler1968unsolved}, Waring's problem on writing integers as sums of $n$th powers \cite{hilbert1909beweis}, and ergodic-theoretic aspects of the Collatz conjecture \cite{lagarias19853x+}. Pisot and Vijayaraghavan proved that this sequence has infinitely many accumulation points \cite{pisot1938repartition}\cite{vijayaraghavan1941fractional}. It is thought to be uniformly distributed, but it has never even been proven dense in the unit interval \cite{finch2003mathematical}. Thus, whether $(3/2)^n$ modulo 1 is uniformly distributed or not is still an unsolved problem.

One straightforward way to understand this sequence is to compute powers of $3/2$ and statistically evaluate the distribution of fractional parts. When generating these powers by multiplication, they quickly result in large numbers that a computer cannot fully represent in floating point. Multiplication is also a computationally expensive operation. To that end, for the problem of efficiently generating powers of $3/2$ modulo 1 we make the following contributions:
\begin{itemize}
  \item We describe a parallelizable, binary addition algorithm of complexity $\mathcal{O}(n^2)$ to compute the fractional part of the sequence $(3/2)^n$ (section \ref{Concept}, appendix \ref{Source}).
  \item We compute the results to $n=10^8$ and conduct statistical testing procedures to assess the strength of the evidence against the uniformity hypothesis for the distribution of the sequence in the unit interval (section \ref{Results}).
\end{itemize}
Section \ref{Background} presents background, section \ref{Concept} offers a description of the proposed parallelization algorithm, section \ref{Results} describes the testing statistical testing procedures employed, section \ref{Related} discusses related work, and section \ref{Conclusion} offers conclusions.

\subsection*{Acknowledgements}
\label{know} J.J.P. Veerman wishes to thank the Department of Mathematics at the University of Crete for its generous hospitality and in particular to Nikos Frantzikinakis for introducing him to this subject.

\section{Background}
\label{Background}
Any real number $x$ can be decomposed into the sum of its parts by letting $\floor*{x}$ denote the integral part and $\{x\} = x - \floor*{x}$ denote the fractional part.

\theoremstyle{definition}
\begin{definition}
A sequence $(x_n)_{n=1}^\infty$ is said to be \textit{uniformly distributed modulo }1 (abbreviated u.d. mod 1) if the proportion of terms of the sequence $\{x_n\}_{n=1}^\infty$ falling in a subinterval $[a,b]$ of $[0,1]$ is proportional to the length of the subinterval \cite{lerma1995distribution}.
That is, if $$\lim_{N\rightarrow \infty} \displaystyle{\frac{\#\{n \leq N: \{x_n\} \in [a,b]\}}{N}} = b-a $$
for every $0 \leq a < b \leq 1.$ Here, $\#(A)$ denotes the cardinality of a set $A$.
\end{definition}

\noindent The next theorem is due to Weyl \cite{weyl1916gleichverteilung}.
\theoremstyle{theorem}
\begin{theorem}
For almost all $x > 1$ and for all $n \in \mathbb{N}$, the sequence $x^n$ is u.d. mod 1.
\end{theorem}
\noindent The exceptional set for which this theorem does not hold has Lebesgue measure 0. A few examples of such sequences in the exceptional set are outlined below, as well as in Finch and Lerma \cite{finch2003mathematical}\cite{lerma1995distribution}:
\begin{enumerate}
  \item All integers $x > 1$, since $x^n \equiv 0$ (mod 1).
  \item Pisot-Vijayaraghavan (or P.V.) numbers, real algebraic integers $x > 1$ whose Galois conjugates lie inside the open unit disc $\{z \in \mathbb{C}: |z| < 1 \}$. Clearly, integers greater than 1 are P.V. numbers. If $x$ is a P.V number, then $\lim_{n\rightarrow \infty} x^n \equiv 0$ (mod 1) geometrically \cite{salem1963algebraic}. Examples are the silver ratio $\delta_s = 1 + \sqrt{2}$ and the golden ratio $\phi = \frac{1 + \sqrt{5}}{2}$.
\end{enumerate}
Remarkably, candidate sequences of the measure 1 set for which $x^n$ is u.d. mod 1 are easy to find but difficult to prove. One such candidate is constructed in Drmota \cite{drmota2006sequences}.

\section{Program Design}
\label{Concept}
\subsection{Binary Representation}
\label{Binary}
A simple relationship exists in binary between the sequences $3^n$ and $(3/2)^n$. The number 3/2 is represented in binary as 1.1. A standard binary addition algorithm to compute $(3/2)^2$ involves shifting the radix point one place to the left and then adding the original binary value to the shifted value. This gives $$1.10 + .11 = 10.01.$$

Similarly, the integer 3 is represented in binary digits as 11. To find $3^2 = 9$, a bit shift operation $(<<)$ can be applied to 11 to obtain 110. This term is then added to the original term, which gives $$110 + 11 = 1001.$$

By continuing this way we obtain ordered sequences of binary numbers representing powers of $(3/2)$ and $3$. Here are the first ten.

\begin{center}
\begin{tabular}{ r c r c }
 1.1 & $(3/2)^1 \quad \quad $ & 11 & $3^1$ \\
10.01 & $(3/2)^2 \quad \quad $ & 1001 & $3^2$ \\
11.011 & $(3/2)^3 \quad \quad $ & 11011 & $3^3$ \\
101.0001 & $(3/2)^4 \quad \quad $ & 1010001 & $3^4$ \\
111.10011 & $(3/2)^5 \quad \quad $ & 11110011 & $3^5$ \\
1011.011011 & $(3/2)^6 \quad \quad $ & 1011011011 & $3^6$ \\
10001.0001011 & $(3/2)^7 \quad \quad $ & 100010001011 & $3^7$ \\
11001.10100001 & $(3/2)^8 \quad \quad $ & 1100110100001 & $3^8$ \\
100110.011100011 & $(3/2)^9 \quad \quad $ & 100110011100011 & $3^9$ \\
111001.1010101001 & $(3/2)^{10} \quad \quad $ & 1110011010101001 & $3^{10}$
\end{tabular}
\end{center}

These sequences produce the same binary digits, the only difference being the inclusion of the radix point in the calculation of $(3/2)^j$ for $1 \leq j \leq n$. A shift of one position in the radix occurs for every iteration because the trailing digit of 1 is always added to 0. Thus, in order to determine $(3/2)^j$ from the calculation of $3^j$, the radix point should be located $j$ positions from the right.

For our purposes, the advantages of generating powers of $3$ instead of powers of $(3/2)$ are twofold. First, a hardware implementation of binary integer addition can be used rather than a software implementation of string addition in order to keep track of the location of the radix point. This improves the speed of program execution. Second, the program can run in parallel by providing powers of 3 as starting values. Large starting values of $3^j$ are calculated using an arbitrary precision package.
\subsection{Algorithm Description}
\label{Algorithm}
This algorithm is an arbitrary precision method specifically tailored to compute the sequence $\{(3/2)^n\}$. A technique called \textit{variable length quantity} is used. That is, we concatenate 63 bits of every 64 bit unsigned long integer to form each value of $3^j$, $1 \leq j \leq n$. The 64th bit serves as a marker to determine when a 64 bit unsigned long integer will overflow. This method is preferable to a standard arbitrary precision package for speed of execution.

We begin by seeding an unsigned long vector called \texttt{value} $\in \mathbb{R}^1$ with the number 3. To generate powers of 3, two arithmetic operations are performed: a shift operation $(<<)$ and an add operation $(+)$. These operations are carried out as described in section \ref{Binary} subject to the following conditions: the shift operation sets a shifted bit when overflow of a 63 bit vector word occurs, and the add operation sets a carry bit when overflow of a 63 bit vector word occurs. A new word is added to the vector when the carry or shift bit, or both, indicates. After this occurs the first time, \texttt{value} becomes an element of $\mathbb{R}^2$. At every iteration, a counter is incremented to track the designated location of the radix point to convert the sequence $3^j$ to $(3/2)^j$.

Next, to determine the distribution of the fractional parts of the sequence $(3/2)^j$, the interval $[0,1)$ is subdivided into $2^r$ bins of equal length. The bin number that corresponds to the fractional part of \texttt{value} for each iteration is determined by selecting the $r$ most significant bits of the fractional part, i.e. the digits in positions $[j, j+r-1]$ from the right. The bin number for each iteration is equal to the decimal representation of the digits $[j, j+r-1]$. For each bin, a counter is incremented every time the bin is hit. Thus, the value of the counter gives the number of elements in the bin.

Below we show a portion of the results for the first 40 iterations using $r=2^{10}$ bins. The fractional part to be binned is underlined. Notice, on the last iteration \texttt{value} spans $\mathbb{R}^2$, but the most significant bits of the fractional part of \texttt{value} spans $\mathbb{R}^1$. The source code provides details for the case where the most significant bits of the fractional part spans $\mathbb{R}^2$ or greater (see appendix \ref{Source}).\\

\begin{scriptsize}
\texttt{\raggedleft
000000000000000000000000000000000000000000000000000000000000001\underline{1} \\
00000000000000000000000000000000000000000000000000000000000010\underline{01} \\
0000000000000000000000000000000000000000000000000000000000011\underline{011} \\
000000000000000000000000000000000000000000000000000000000101\underline{0001} \\
00000000000000000000000000000000000000000000000000000000111\underline{10011} \\
0000000000000000000000000000000000000000000000000000001011\underline{011001} \\
000000000000000000000000000000000000000000000000000010001\underline{0001011} \\
00000000000000000000000000000000000000000000000000011001\underline{10100001} \\
0000000000000000000000000000000000000000000000000100110\underline{011100011} \\
000000000000000000000000000000000000000000000000111001\underline{1010101001} \\
00000000000000000000000000000000000000000000001010110\underline{0111111101}1 \\
0000000000000000000000000000000000000000000010000001\underline{1011111100}01 \\
\textbf{\vdots} \quad \quad \quad \quad \quad \quad \quad \quad \quad \quad \quad \quad \quad \quad \quad \quad \\
0011100000111101100100010\underline{1110000101}11000010111111111100000001011 \\
...000000000000000000000000000000000001  \ 001010001011100010110100\underline{0101001000}101001000111111110100000100001 \\
\color{white}hey!}
\end{scriptsize}
\noindent The following table displays the bin numbers for this example using $2^{10}$ bins.

\begin{table}[H]
\centering
\label{table:1}
\captionsetup{justification=centering}
\caption{}
\begin{tabular}{ | c | c | c | }
\hline
Iteration & Most significant fractional digits & Corresponding bin number \\ \hline
1 & 1000000000& 512 \\ \hline
2 & 0100000000& 256 \\ \hline
3 & 0110000000& 384 \\ \hline
4 & 0001000000& 064 \\ \hline
5 & 1001100000& 608 \\ \hline
6 & 0110010000& 400 \\ \hline
7 & 0001011000& 088 \\ \hline
8 & 1010000100& 644 \\ \hline
9 & 0111000110& 454 \\ \hline
10 & 1010101001& 681 \\ \hline
11 & 0111111101& 509 \\ \hline
12 & 1011111100& 764 \\ \hline
\vdots & \vdots& \vdots \\ \hline
39 & 1010101001& 901 \\ \hline
40 & 0101001000& 328 \\ \hline
\end{tabular}
\end{table}

\section{Results}
\label{Results}
We now conduct a numerical experiment for large $n$, and implement a statistical testing procedure to assess the strength of the evidence for the uniformity of the distribution $(3/2)^n$. In particular, we are interested in assigning $n=10^8$ values in $[0,1]$ to $r=2^{25}$ bins that partition the interval, to determine whether the values are uniformly distributed across the different bins.  The number of ways in which $n$ objects can be assigned to $r$ bins ($r\leq n$) is given by the multinomial coefficient $n!/\lrp{x_1!x_2!\cdots x_r!}$.

Hence, assuming that all bins are equally probable (i.e., the numbers are uniformly distributed within the unit interval) and letting $X_i$ ($1\leq i \leq r$) denote the random variable that counts the number of values that fall in bin $i$, the random vector $\boldsymbol{X}=(X_1,\ldots,X_r)$, jointly follows a multinomial distribution with probability mass function (pmf) given by
\begin{eqnarray}
\Pr(X_1=x_1,\ldots,X_r=x_r \given  H_0)&=& \frac{n!}{x_1!x_2!\cdots x_r!} \left(\frac{1}{r}\right)^{x_1+\ldots+x_r}\nonumber\\
&=&{n\choose x_1,\ldots,x_r}\left(\frac{1}{r}\right)^n.\label{eq:H0}
\end{eqnarray}
The conditioning on $H_0$ is used to represent the distribution of $\boldsymbol{X}$ under the null hypothesis of  uniformity. In general , letting $\bpi=(\pi_1,\pi_2,\ldots,\pi_r)$ the multinomial pmf is given by
\begin{eqnarray}
f(\bx \given \bpi)\;=\; \Pr(X_1=x_1,\ldots,X_r=x_r \given \bpi)&=& \frac{n!}{x_1!x_2!\cdots x_r!} \pi_1^{x_1}\cdots\pi_r^{x_r}.\label{eq:H1}
\end{eqnarray}

Thus, testing for uniformity amounts to deciding if the vector $$\bpi=\lrp{\pi_1,\pi_2,\ldots,\pi_r}\text{ is equal to }\bpi_0=\lrp{1/r,1/r,\ldots,1/r},$$ based on the observed vector of counts for the $r$ bins. This statistical inference problem can be addressed from different viewpoints depending on the branch of statistics considered, each providing fundamentally different inferential statements.  From the frequentist perspective, the strength of the evidence \emph{against} $H_0: \bpi=\bpi_0$ can be assessed through a testing procedure that yields a falsifiability statement.  Conversely, the Bayesian inferential setting allows one to compare the evidence present in the data for both $H_0$ and $H_1: \bpi\not=\bpi_0$ through a probabilistic statement.  Throughout the remainder of the section we describe each approach and provide the results obtained for our data in each case.

\subsection{Frequentist Testing}
Defining $\hat{\bpi}=\lrb{X_1/n, X_2/n,\ldots, X_r/n}$, the test statistic $T(\hat{\bpi},\bpi_0)$ is a random variable that characterizes the divergence between the multinomial distributions characterized by the densities \eqref{eq:H0} and \eqref{eq:H1}.  Let the vector of observed counts for the $r$ bins be given by $\bx=(x_1, x_2, \ldots, x_r)$, and denote by $\hat{\bp}=(x_1/n, x_2/n, \ldots, x_r/n)$ the vector of estimated values for the multinomial probabilities.  Furthermore, let $\mcal{R}(\alpha)$ denote the $\alpha$-level rejection region for the test.  Then the testing procedure is characterized by the decision rule
\bean
\text{If }\tau\in\mcal{R}(\alpha) \hfill & \longrightarrow & \text{reject }H_0\\
\text{If }\tau\not\in\mcal{R}(\alpha) \hfill & \longrightarrow & \text{do not reject }H_0,
\eean
where $\tau=T(\hat{\bp},\bpi_0)$ (i.e., the value of the test statistic evaluated at $\hat{\bp}$).  The rejection region $\mcal{R}(\alpha)$ depends on the distribution of the test statistic $T(\hat{\bpi},\bpi_0)$ and the level $\alpha$, often chosen to be in the set $\lrb{0.01,0.05,0.1}$. 

The traditional asymptotic results for the distribution of $T(\hat{\bpi},\bpi_0)$ are obtained for the case with $r/n\rightarrow 0$ as $n\rightarrow\infty$ \cite{read2012goodness,Bahadur08}.  However, our problem of interest includes an exceedingly large number of bins, hence, a more suitable assumption for the asymptotic relationship between $n$ and $r$ is $r/n\rightarrow \nu$ with $\nu\in(0,\infty)$. In particular, \cite{Holst1972} shows that some traditional test statistics for the $r/n\rightarrow 0$ case, such as the likelihood ratio test, break down if $r/n\rightarrow \nu$; however, \cite{Holst1972} demonstrates that the popular $\chi$-squared statistic, given by
\begin{eqnarray}
T(\hat{\bpi},\bpi_0)&=& \lrp{\frac{r}{n}\sum_{i=1}^r X_i^2}-n \;=\; rn\lrb{\lrp{\sum_{i=1}^r \lrp{\frac{X_i}{n}}^2}-\frac{1}{r}},
\end{eqnarray}
is robust to this type of asymptotic behavior.  Assuming that $H_0$ holds, $T(\hat{\bpi},\bpi_0)$ is asymptotically distributed chi-squared with $r-1$ degrees of freedom \cite{pearson1900x} (denoted by $\chi_{r-1}^2$) as $r/n\rightarrow\nu$ while $n$ and $r$ approach $\infty$. Hence, for this test statistic the $\alpha$-level rejection region is given by $\mcal{R}(\alpha)=[\chi_{r-1,1-\alpha},\infty)$, where $\chi_{r-1,1-\alpha}$ denotes the $\alpha$-level critical value for the $\chi_{r-1}^2$ distribution (e.g., see Chapter 8.2 in \cite{casella2002statistical}). 
Hence, the test statistic evaluated at $\hat{\bp}$ corresponds to
\bea
\tau \;=\; T(\hat{\bp},\bpi_0) &=& \lrp{\frac{r}{n}\sum_{i=1}^r x_i^2}-n,\label{eq:tau}
\eea
implying that if $\hat{p}_i=\frac{x_i}{n}=\frac{1}{r}$ for all $i=1,\ldots,r$, then $\tau=0$, which explains why the rejection region only encompasses the upper tail of the distribution.

Figure \ref{fig:chisqtest} below depicts the probability density functions obtained for the test statistic under the null hypothesis $H_0$ for $r=2^{25}$ and $n=10^8$.  The region shaded in grey represents $P\lrp{T(\hat{\bpi},\bpi_0) > \tau}$.
As mentioned above, the $\alpha$-level hypothesis test conducted rejects $H_0$ whenever $T(\hat{\bpi},\bpi_0)\in\mcal{R}(\alpha)$, or equivalently whenever $P\lrp{T(\hat{\bpi},\bpi_0) > \tau} < \alpha$, which is clearly not the case for for $\alpha\in\lrb{0.01, 0.05,0.1}$ (see Figure \ref{fig:chisqtest}). Thus, from the frequentist standpoint there is no evidence in the data for the test to reject the null hypothesis in favor of $H_1$.

\begin{figure}[htbp]
\begin{center}
\includegraphics[scale=0.5]{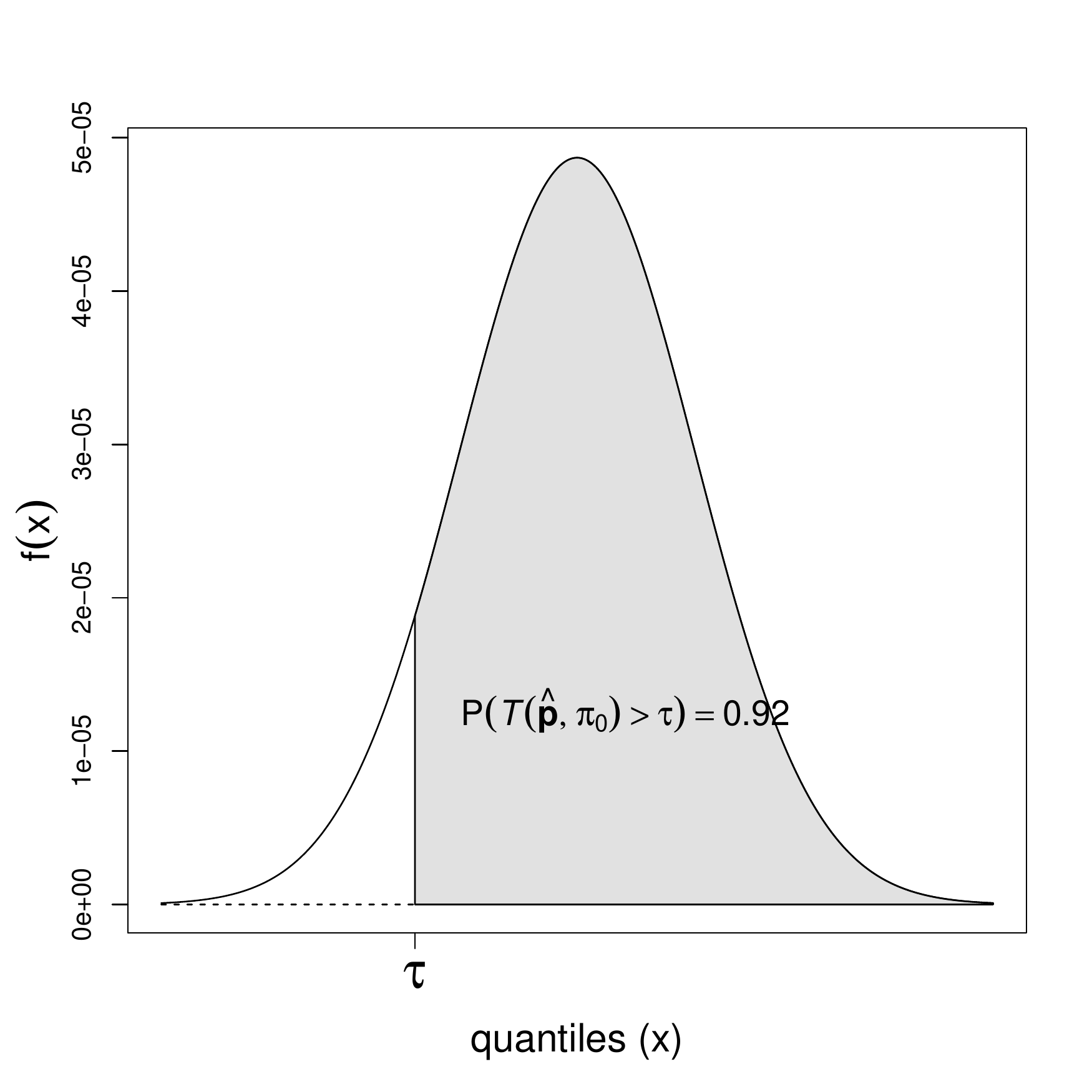}
\caption{Probability density function for $T(\hat{\bp},\bpi_0)$ under the null hypothesis of uniformity with $n=10^8$ and $r=2^{25}$.}
\label{fig:chisqtest}
\end{center}
\end{figure}

We tried analogous tests for the pairs $(n=(10^8/2^3), r=2^{22})$, $(n=(10^8/2^5), r=2^{20})$,  and $(n=(10^8/2^7), r=2^{18})$, with $P\lrp{T(\hat{\bpi},\bpi_0) > \tau}$ equal to 0.45, 0.94, and 0.01, respectively.  Note that the three tests with larger $n$ and $r$ values, lead to probabilities that are consistent with $H_0$.  While the test for $(n=(10^8/2^7), r=2^{18})$ has a small probability under $H_0$, this result should be interpreted in contrast to all other possible multinomial weights. In the next section we will revisit this result from the Bayesian standpoint, which allows us to test for the strength of $H_0$ among all other possibilities in the $r$-dimensional simplex, as opposed to simply providing an argument of the falsibiability of $H_0$.

\subsection{Bayesian Inference}

While the frequentist testing framework assumes that the vector of multinomial probabilities $\bpi$ is fixed but unknown, in the Bayesian context $\bpi$ is considered a random vector. The Bayesian testing machinery uses Bayes rule to calculate the conditional probability for the hypotheses of interest given the observed data, denoted by $P(H_0\given \bx)$ and $P(H_1\given \bx)$, respectively, and uses these probabilities as the measure of evidence for or against them.  These probabilities are referred to in the statistical literature as \emph{posterior probabilities}.  The difficulty in calculating these probabilities stems from choosing \emph{suitable} prior probability distributions over the $r$-dimensional simplex $\mcal{S}_r$ (for the vector of multinomial probabilities $\bpi$) and over the space of hypotheses $\lrb{H_0,H_1}$.   What is a \emph{suitable} prior probability distribution has been a topic of intense debate in the statistical literature, but the consensus is that in the absence of prior information, these statements should avoid biasing the outcome of the tests in any direction, and should thus be as uninformative as possible to maximize the ability of the information contained in the data to drive the conclusions.

For our particular problem selecting a prior distribution on the space of hypotheses is straightforward. This space solely contains two elements, therefore, only $\psi_0=P(H_0)$ (the prior probability that  $H_0$ holds) has to be specified.  By construction, conditional on $H_0$ the vectors $\bpi$ and $\bpi_0$ are equal with probability one. Furthermore, let $g$ denote a conditional  density function with respect to the Lebesgue measure on $\lrb{\bpi\not=\bpi_0}$, with $g(\bpi_0)=0$.  Using Bayes theorem and the definitions above, the posterior probability for $H_0$ is
\bea
P(H_0\given \bx)&=&\frac{\psi_0 \, f(\bx\given \bpi_0) }{\psi_0 \, f(\bx\given \bpi_0) + (1-\psi_0) \int_{\mcal{S}_r}f(\bx\given \bpi)g(\bpi)d\,\bpi}\nonumber\\
&=& \lrp{1+\frac{(1-\psi_0)}{\psi_0} \frac{1}{B_g(\bx)}}^{-1},\label{eq:post}
\eea
where $B_g(\bx)=f(\bx\given \bpi_0)\big/ \int_{\mcal{S}_r}f(\bx\given \bpi)g(\bpi)d\,\bpi$ is the Bayes factor, which corresponds to a likelihood ratio for $H_0$ relative to $H_1$ (where the likelihood under $H_1$ is averaged over all possible $\bpi\in\mcal{S}$ according to $g$).

In the absence of prior information the usual choice for $\psi_0$ is often $\psi=0.5$, since it sets the ratio in $(1-\psi_0)/\psi_0$ in \eqref{eq:post}  to 1, making
\bean
P(H_0\given \bx) &=& \frac{B_g(\bx)}{1+ B_g(\bx)},
\eean
which is a function driven by the observed data and the function $g$.  The choice of $g$ strongly affects the outcome of the testing problem, with its influence manifesting through $B_g(\bx)$.  Given that $P(H_0\given \bx)$ is a monotone nondecreasing function of $B_g(\bx)$, instead of choosing a particular function $g$, in \cite{Delampady1990} the authors propose to obtain a lower bound $\underline{B}_G(\bx)=\inf_{g\in G}B_g(\bx)$ over entire classes of densities $G$.  This in turns determines a conservative lower bound for $P(H_0\given \bx)$.  Here, we follow the guidance provided in \cite{Delampady1990}, and restrict to the Dirichlet distribution, characterized by the class of density functions of the form
$$G=\lrb{ g_{c}(\bpi)=\frac{\Gamma(c)}{\lrp{\Gamma( c/r)}^r} \prod_{i=1}^r\pi_{i}^{(c/r)-1}\text{, for some }c>0\text{ and }\bpi\in\mcal{S}_r},$$
whose elements $g_{c}(\bpi)$ are conjugate to the likelihood $f(\bx\given \bpi)$.  These density functions are centered around $\bpi_0$ (i.e., $E_{g_c}[\bpi]=\bpi_0$), and spread their mass around the mean in a way that prevents biases towards particular values, lending some measure of ``objectivity'' to the testing procedure. The lower bound for the Bayes factor in the class $G$ is given by
\bean
\underline{B}_G(\bx) &=&\inf_{c>0}\frac{\Gamma(c+n)\lrp{\Gamma(c/r)}^r \lrp{1/r}^{n}}{\Gamma(c)\prod_{i=1}^r\Gamma(x_i+c/r)},
\eean
which approaches $\underline{B}^{*}(\bx)=\inf_{0<\rho<1} \rho^{1-r}\exp{\lrb{-\frac{1}{2}(1-\rho^2)\tau}}$ as $n\rightarrow\infty$, with $\tau$ defined as in \eqref{eq:tau}.

Applying the results from \cite{Delampady1990} described above to the data generated with $n=10^8$ and $r=2^{25}$, the value of the asymptotic lower bound $\underline{B}^{*}(\bx)=1$.  A Bayes factor of 1 represents indifference between two competing hypotheses.  That being said, the fact that the Bayes factor is bounded below at 1 indicates that, at its worst, the evidence for $\bpi_0$ is as strong as the evidence in favor all other possible values of $\bpi\in \lrb{\mcal{S}_r \setminus \lrb{\bpi_0}}$ combined.   Thus, choosing $\psi_0=0.5$, we have that
\bean
P(H_0\given \bx) &\geq& \frac{B^*(\bx)}{1+ B^*(\bx)} \; =\; \frac{1}{2},
\eean
implying that the posterior probability for $\lrb{\bpi=\bpi_0}$ is at least as large as the posterior probabilities for all other possible values in $\mcal{S}_r$ combined; an overwhelming amount of evidence in favor of $H_0$.

Revisiting the alternative experiments considered, we have that, excluding the pair $n=(10^8/2^7)$ and $r=2^{18}$, all other $(n,r)$ combinations lead lower bound for $P(H_0\given \bx)$ of $1/2$, which again provides strong evidence in favor of $H_0$.  For the pair $n=(10^8/2^7)$ and $r=2^{18}$, we find a lower bound of 8\%.  Implying that although the evidence is not overwhelming in this case, there is some evidence for the null hypothesis, especially in light of the fact that the space of all other possibilities is uncountable.

\section{Related Work}
\label{Related}
To our knowledge no research has been published to evaluate numerically the distribution of the sequence $(3/2)^n$ modulo 1. However, the distribution of this sequence has been studied in relation to another topic: Mahler's Z-numbers. A Z-number is a real number $\xi$ such that $\{\xi(3/2)^n\} \leq 1/2$ for all $n \in \mathbb{N}$. Mahler conjectured in 1968 that no Z-numbers exist \cite{mahler1968unsolved}.

Recently, Akiyama proved that for coprime integers $p>q> 1$, where $p > q^2$, there exists $\xi > 0$ such that $\{\xi(p/q)^n\}$ stays in a Cantor set for all $n \in \mathbb{N}$ \cite{akiyama2008mahler}. Furthermore, Dubickas also proved that for coprime $p>q> 1$, there exists $\xi > 0$ such that $\{\xi(p/q)^n\}$ lies in a short interval of $[0,1]$. In particular, his work implies that $\Vert \{\xi (3/2)^{2n}\} \Vert < 14/45$ \cite{dubickas2010powers}.

Waring's problem is the conjecture that every natural number can be expressed as the sum of $k$ $n$th powers of non-negative integers \cite{finch2003mathematical}. Waring's problem can be solved if the inequality $$\{(3/2)^n\} \leq 1 - (3/4)^n$$ holds. Kubina and Wunderlich demonstrated computationally that this inequality is met for $2 \leq n \leq 471,600,000$ \cite{kubina1990extending}. No counterexample to this inequality has been found, and it is believed that this inequality can be extended to $$(3/4)^n < \{(3/2)^n\} < 1 - (3/4)^n$$ for $n > 7$, though this statement has never been proven \cite{finch2003mathematical}.

As Vijayaraghavan \cite{vijayaraghavan1941fractional} remarked in 1941, it is not known whether $$\limsup\limits_{n\rightarrow \infty} \{(3/2)^n\} - \liminf\limits_{n\rightarrow \infty} \{(3/2)^n\} > 1/2,$$ and this is still an open problem, but in 1995, Flatto, Lagarias, and Pollington \cite{flatto1995range} showed that $$\limsup\limits_{n\rightarrow \infty} \{(3/2)^n\} - \liminf\limits_{n\rightarrow \infty} \{(3/2)^n\} \geq 1/3.$$

\section{Conclusions and Further Work}
\label{Conclusion}
The statistical tests conducted for the numerical experiment strongly support the hypothesis that the sequence $(3/2)^n$ modulo 1 is uniformly distributed when larger sample sizes are considered, both from the frequentist and from the Bayesian standpoint.  For a smaller dataset tested, the evidence is not as conclusive, but even in this case the uniformity hypothesis is not unlikely when considering the outcome of the Bayesian testing procedure.

For the numerical experiment described in Section \ref{Results}, the number of iterations was taken to $n= 10^8$, but it could reasonably be taken to $n=10^9$ or $10^{10}$ with the aid of a supercomputer. The program was parallelized to run four instances total on two computers with Intel(R) Core(TM) i7-2600 CPU @ 3.40GHz, 3401 Mhz, 4 Core(s), 8 Logical Processor(s). The runtime of the program for these computers was $11$ days, with a more precise runtime given by $t = (2.378 \num{1e-10})n^2$, for $n$ denoting the number of iterations and $t$ denoting time (in seconds). Starting values for parallelization are noted in Figure \ref{speed} below.

\begin{figure}[H]
\centering
\includegraphics[width=\textwidth]{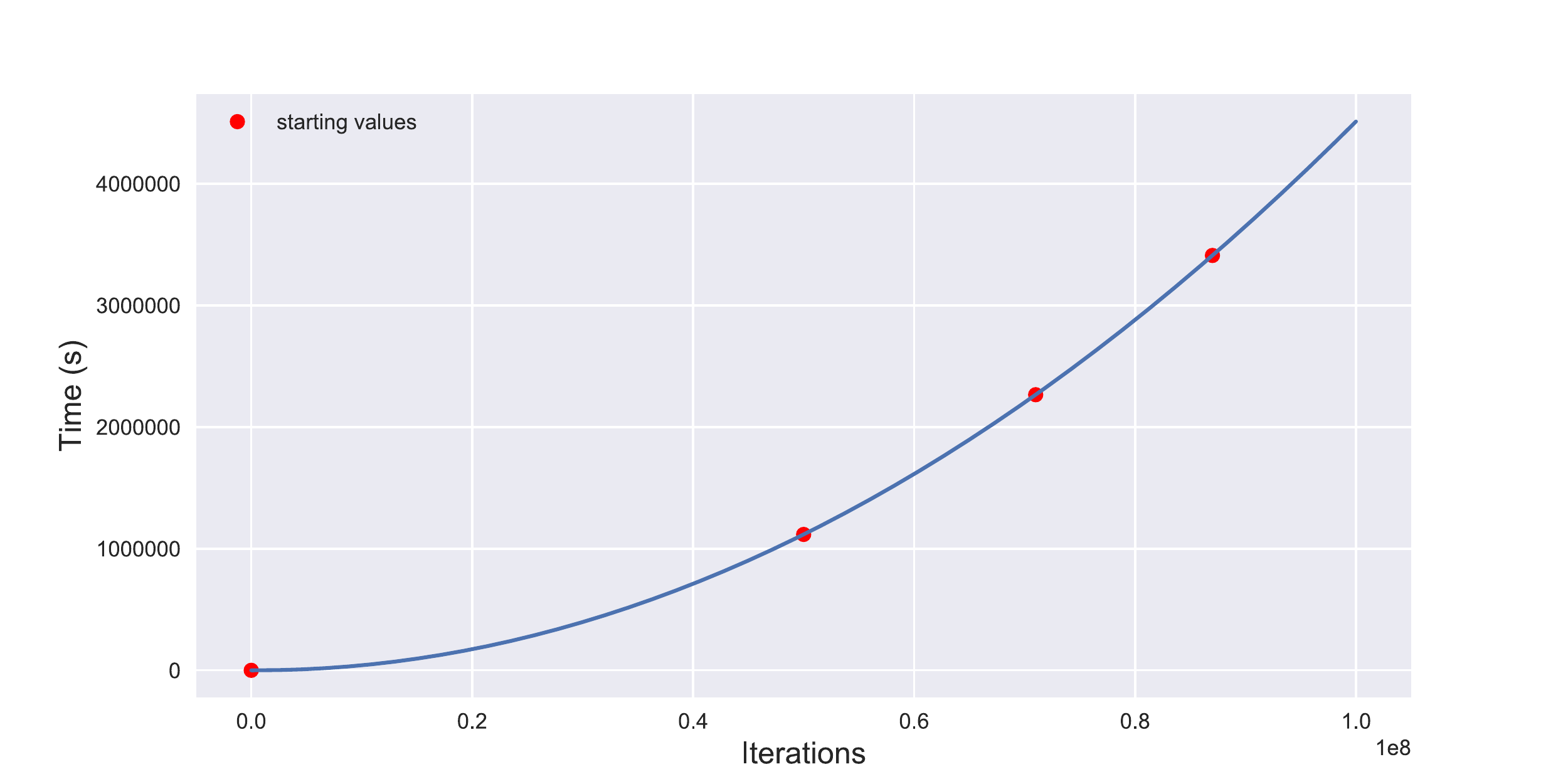}
\caption{Program runtime (extrapolated)}
\label{speed}
\end{figure}

\bibliography{references}

\begin{thebibliography}{10}

\bibitem{akiyama2008mahler}
{\sc Akiyama, S.}
\newblock Mahler’s z-number and 3/2 number systems.
\newblock {\em Uniform Distribution Theory 3}, 2 (2008), 91--99.

\bibitem{casella2002statistical}
{\sc Casella, G., and Berger, R.~L.}
\newblock {\em Statistical Inference}, vol.~2.
\newblock Duxbury Pacific Grove, CA, 2002.

\bibitem{Delampady1990}
{\sc Delampady, M., and Berger, J.~O.}
\newblock {Lower bounds on bayes factors for multinomial distributions, with
  application to chi-squared tests of fit}.
\newblock {\em The Annals of Statistics 18}, 3 (1990), 1295--1316.

\bibitem{drmota2006sequences}
{\sc Drmota, M., and Tichy, R.~F.}
\newblock {\em Sequences, Discrepancies and Applications}.
\newblock Springer, 2006.

\bibitem{dubickas2010powers}
{\sc Dubickas, A.}
\newblock Powers of rational numbers modulo 1 lying in short intervals.
\newblock {\em Results in Mathematics 57}, 1 (2010), 23--31.

\bibitem{finch2003mathematical}
{\sc Finch, S.~R.}
\newblock {\em Mathematical Constants}, vol.~93.
\newblock Cambridge university press, 2003.

\bibitem{flatto1995range}
{\sc Flatto, L., Lagarias, J.~C., and Pollington, A.~D.}
\newblock On the range of fractional parts $\{\xi(p/q)^n\}$.
\newblock {\em Acta Arith 70}, 2 (1995), 125--147.

\bibitem{Bahadur08}
{\sc Harremo{\"e}s, P., and Vajda, I.}
\newblock On the bahadur-efficient testing of uniformity by means of entropy.
\newblock {\em IEEE Transactions on Information Theory 54\/} (2008), 321--331.

\bibitem{hilbert1909beweis}
{\sc Hilbert, D.}
\newblock Beweis f{\"u}r die \uppercase{D}arstellbarkeit der ganzen
  \uppercase{Z}ahlen durch eine feste \uppercase{A}nzahln ter
  \uppercase{P}otenzen (\uppercase{W}aringsches problem).
\newblock {\em Mathematische Annalen 67}, 3 (1909), 281--300.

\bibitem{Holst1972}
{\sc Holst, L.}
\newblock {Asymptotic normality and efficiency for certain goodness-of-fit
  tests}.
\newblock {\em Biometrika 59}, 1 (1972), 137--145.

\bibitem{koksma}
{\sc Koksma, J.~F.}
\newblock \uppercase{E}in mengentheoretischer \uppercase{S}atz {\"u}ber die
  \uppercase{G}leichverteilung modulo \uppercase{E}ins.
\newblock {\em Compositio Mathematica 2\/} (1935), 250--258.

\bibitem{kubina1990extending}
{\sc Kubina, J.~M., and Wunderlich, M.~C.}
\newblock Extending \uppercase{W}aring’s conjecture to 471,600,000.
\newblock {\em Mathematics of Computation 55}, 192 (1990), 815--820.

\bibitem{lagarias19853x+}
{\sc Lagarias, J.~C.}
\newblock The 3x+ 1 problem and its generalizations.
\newblock {\em The American Mathematical Monthly 92}, 1 (1985), 3--23.

\bibitem{lerma1995distribution}
{\sc Lerma, M.~A.}
\newblock Distribution of powers modulo 1 and related topics.

\bibitem{mahler1968unsolved}
{\sc Mahler, K.}
\newblock An unsolved problem on the powers of 3/2.
\newblock {\em Journal of the Australian Mathematical Society 8}, 2 (1968),
  313--321.

\bibitem{pearson1900x}
{\sc Pearson, K.}
\newblock On the criterion that a given system of deviations from the probable
  in the case of a correlated system of variables is such that it can be
  reasonably supposed to have arisen from random sampling.
\newblock {\em The London, Edinburgh, and Dublin Philosophical Magazine and
  Journal of Science 50}, 302 (1900), 157--175.

\bibitem{pisot1938repartition}
{\sc Pisot, C.}
\newblock La r{\'e}partition modulo 1 et les nombres alg{\'e}briques.
\newblock {\em Annali della Scuola Normale Superiore di Pisa-Classe di Scienze
  7}, 3-4 (1938), 205--248.

\bibitem{read2012goodness}
{\sc Read, T.~R., and Cressie, N.~A.}
\newblock {\em Goodness-of-fit Statistics for Discrete Multivariate Data}.
\newblock Springer Science \& Business Media, 2012.

\bibitem{salem1963algebraic}
{\sc Salem, R.}
\newblock {\em Algebraic Numbers and Fourier Analysis}.
\newblock Wadsworth Company, 1963.

\bibitem{vijayaraghavan1941fractional}
{\sc Vijayaraghavan, T.}
\newblock On the fractional parts of the powers of a number (ii).
\newblock In {\em Mathematical Proceedings of the Cambridge Philosophical
  Society\/} (1941), vol.~37, Cambridge University Press, pp.~349--357.

\end{thebibliography}
\bibliographystyle{acm}

\appendix

\section{Source Code}
\label{Source}
All source code can be found at www.github.com/paulaneeley.

\end{document}